\documentclass[12pt]{article}
\usepackage{amsmath}
\usepackage{amsfonts}
\usepackage{amssymb}
\usepackage{amsthm}
\usepackage{hyperref}

\usepackage{graphicx}

\def\bee{\begin{equation}}
\def\eee{\end{equation}}
\def \BB{${\mathcal{B}_2}$}
\def \BBB{{\mathcal{B}_2}}
\def\Li{{\rm {L}i}}

\begin{document}

\thispagestyle{empty}
\centerline{}
\bigskip
\bigskip
\bigskip
\bigskip
\bigskip
\bigskip
\centerline{\Large\bf Remark on the irrationality of the  Brun's constant}
\bigskip
\bigskip

\begin{center}
{\large \sl Marek Wolf}\\*[5mm]
e-mail:mwolf@ift.uni.wroc.pl\\
\bigskip
\end{center}

\bigskip
\bigskip

\begin{center}
{\bf Abstract}\\
\bigskip
\begin{minipage}{12.8cm}
We have calculated  numerically geometrical means of the denominators of
the continued fraction approximations to the   Brun constant
\BB.  We get values close to the
Khinchin's constant. Next we calculated the $n$-th square roots of the
denominators of the $n$-th convergents of these continued fractions obtaining
values close to the Khinchin-L{`e}vy constant. These two results suggests
that \BB  ~ is irrational, supporting  the common believe that there is an infinity of twins.
\end{minipage}
\end{center}

\bigskip\bigskip

\bibliographystyle{abbrv}

Very well known open problem in number theory  is the question whether there exist
infinitely many twin primes $p, p+2$.  In  1919 Brun \cite{Brun} has shown that the sum of the
reciprocals of all twin primes is finite:
\bee
\BBB=\left( {1 \over 3} + {1\over 5}\right) + \left( {1 \over 5} + {1\over
7}\right) + \left( {1 \over 11} + {1\over 13}\right) + \ldots < \infty,
\label{def-B2}
\eee
thus leaving the problem  not decided. Sometimes 5 is included in (\ref{def-B2})
only once, but here we  will adopt the above
convention. 
The sum (\ref{def-B2}) is called the Brun constant \cite{Wrench1974}.

Let $\pi_2(x)$ denote the number of twin primes smaller than $x$. Then
the conjecture B of Hardy and Littlewood \cite{Hardy_and_Littlewood} on
the number of prime pairs $p, p+d$ applied to the case $d=2$ gives, that
\begin{equation}
\pi_2(x) \sim C_2\Li_2(x) \equiv C_2 \int_2^x \frac{u}{\log^2(u)} du,
\label{conj}
\end{equation}
where $C_2$ is called ``twin constant'' and is defined by the following
infinite product:
\begin{equation}
C_2 \equiv 2 \prod_{p > 2} \biggl( 1 - \frac{1}{(p - 1)^2}\biggr) =
1.32032363169\ldots
\label{stalac2}
\end{equation}
There is a large evidence  both analytical and experimental in favor of (\ref{conj}).
Besides the original circle method used by Hardy and Littlewood
\cite{Hardy_and_Littlewood}  there appeared
the paper \cite{Rubinstein} where another heuristic arguments were presented.  
In May 2004, in a preprint publication \cite{Arenstorf} Arenstorf attempted to prove
that there are infinitely many twins. Arenstorf  tried to continue analytically
to  $\Re (s)=1$ the difference:
\bee
T(s)-C_2/(s-1)
\eee
where the function
\bee
T(s)= \sum_{n>3}\Lambda(n-1)\Lambda(n+1)n^{-s} \quad (\Re  (s) > 1).
\label{zeta_twins}
\eee
However shortly after an error in the proof was
pointed out by Tenenbaum \cite{Tenenbaum}.  For  recent progress in the direction of
the proof of the infinite number of twins see \cite{Koreevar}.
Because  there is no doubt that twins prime conjecture is true the Brun's constant
should be irrational.

The series (\ref{def-B2}) is very slowly convergent and there is a method based on
the (\ref{conj}) to extrapolate finite size approximations
\bee
\BBB(x)=\sum_{q, q+2 ~{ both ~prime} \atop q \leq x} \Big(\frac{1}{q} + \frac{1}{q+2}\Big)
\eee
to infinity \cite{Brent}:
\bee
\BBB^{\!\!(\infty)} (x)=\BBB(x) + \frac{2C_2}{\log(x)}
\eee
In this way from the straight sieving of primes  up to $x=3\times 10^{15}$ T. Nicely
\cite{Nicely_Brun} gives
\bee
\BBB ^{\!\!(\infty)} (3\times 10^{15}) = 1.90216 05823 \pm 0.00000 00008
\eee
while P. Sebah \cite{Sebah-Brun} from computer search up to $10^{16}$ gives
\bee
\BBB^{\!\!(\infty)} (10^{16}) = 1.90216 05831 04.
\label{Sebah}
\eee


If there is an infinity of twins, as the formula (\ref{conj}) asserts,
then the Brun's constant should be an irrational number.
Vice versa if the Brun's constant is irrational then there is an
infinity of twins.

There exists a method based on the continued fraction expansion which
allows  to detect whether a given number $r$ can be irrational or not.   Let
\bee
r = [a_0; a_1, a_2, a_3, \ldots]=
a_0+\cfrac{1}{a_1 + \cfrac{1}{a_2 + \cfrac{1}{a_3+\ddots}}}
\eee
be the continued fraction expansion of the real number $r$, where $a_0$ is an
integer and  all $a_k,~k=1, 2, \ldots$ are
positive integers. Khinchin has proved that
\bee
\lim_{n\rightarrow \infty} \big(a_1 a_2 \ldots a_n\big)^{\frac{1}{n}}=
\prod_{m=1}^\infty {\left\{ 1+{1\over m(m+2)}\right\}}^{\log_2 m} \equiv K_0 \approx 2.685452001\ldots
\eee
is a constant for almost all real $r$, see e.g. \cite[\S 1.8]{Finch}.
The exceptions are {\it rational numbers},
quadratic irrationals and some irrational numbers too, like for example the
Euler constant $e=2.7182818285\ldots$, but this set of exceptions is of the
Lebesgue measure zero. The constant $K_0$ is called the Khinchin constant.
If the quantities
\bee
K(n)=\big(a_1 a_2 \ldots a_n\big)^{\frac{1}{n}}
\eee
for   $r$  given with  accuracy of some number of digits are close to $K_0$
we can regard it as a hint that  $r$ is irrational.

For the numerical value of \BB ~ given by  (\ref{Sebah}) we get continued fraction
containing 23 terms:
\[
\BBB\approx 1.902160583104=[1; 1, 9, 4, 1, 1, 8, 3, 4, 4, 2, 2, 2, 1, 35, 1, 1, 1, 2, 4, 4, 1, 2]
\]
\bee
=1+\cfrac{1}{1+\cfrac{1}{9+\cfrac{1}{4+\cfrac{1}{1+\cfrac{1}{1+\cfrac{1}{1+\cfrac{1}{8+\cfrac{1}{3+\cfrac{1}{4+\cfrac{1}{4+\ddots}}}}}}}}}}
\label{cfr-Brun}
\eee
We have calculated the geometrical means $K(n)$  for the  consecutive
truncations of the continued fraction (\ref{cfr-Brun})  for $n=7, 8, \ldots 23$.
The results are presented in the second column of Table 1 and in Fig.1 for $n=3, 4,
\ldots, 23$ and $K(n)$
are fluctuating around $K_0$, suggesting \BB  \  is indeed  irrational.

\newpage

\vskip 0.4cm
\begin{center}
{\sf TABLE {\bf I}}\\
\bigskip
\begin{tabular}{||c|c|c||} \hline
$n$ & $ K(n) $ & $ L(n) $ \\ \hline
$7$ & 2.5697965868506505913  & 3.038522491494629559394  \\ \hline
$8$ & 2.6272534028385753462    & 3.049597387327883630584 \\ \hline
$9$ & 2.7689921001973389000    & 3.185356988368901235392 \\ \hline
$10$ & 2.8844991406148167646    & 3.287467122665676043642 \\ \hline
$11$ & 2.7807783406318186810    & 3.163215603753400545370 \\ \hline
$12$ & 2.6986960551881232673    & 3.090263509588119664954 \\ \hline
$13$ & 2.6321480259049849216    & 3.026751930527273995360 \\ \hline
$14$ & 2.4433109823349210157  & 2.854823973874123136328 \\ \hline
$15$ & 2.9549936641069114440    & 3.419389588999534265980 \\ \hline
$16$ & 2.7490720167650636040    & 3.156107031770304366266 \\ \hline
$17$ & 2.5806968087738011680    & 3.064762600086287594909 \\ \hline
$18$ & 2.4407155166256854482 &  2.939414938043217059705 \\ \hline
$19$ & 2.4138613420275263905 &  2.923369271618110487047 \\ \hline
$20$ & 2.4788881641188601981  &  2.986079858472829061265 \\ \hline
$21$ & 2.5389087172991914976    & 3.038476735121596843751 \\ \hline
$22$ & 2.4287240750956647922  & 2.911134637211735156321 \\ \hline
$23$ & 2.4073773469514141247  & 2.848836285543623055400 \\ \hline
\end{tabular} \\
\end{center}
\vskip 0.4cm

Let the rational $p_n/q_n$ be the $n$-th  partial convergent of the continued fraction:
\bee
\frac{p_n}{q_n}=[a_0; a_1, a_2, a_3, \ldots, a_n].
\eee
For almost all real numbers $r$ the denominators of the  finite continued fraction
approximations fulfill:
\bee
\lim_{n \rightarrow \infty} \left(q_n(r)\right)^{1/n} = e^{\pi^2/12\ln2} \equiv L_0 =  3.275822918721811\ldots
\eee
where $L_0$ is called the  Khinchin---L{\`e}vy's constant \cite[\S 1.8]{Finch}.
Again the set of exceptions to
the above limit is  of the Lebesgue measure zero and it includes {\it rational numbers},
quadratic irrational etc.

From (\ref{cfr-Brun}) we get the following sequence  of  convergents $p_n(\BBB)/q_n(\BBB)$:
$$
\bigg( \frac{1} {1} , ~~  \frac{2} {1} , ~~ \frac{19} {10} , ~~  \frac{78} {41} , ~~  \frac{97} {51} , ~~  \frac{175} {92} , ~~  \frac{1497} {787} , ~~  \frac{4666} {2453} , ~~  \frac{20161} {10599} , ~~  \frac{85310} {44849} , ~~  \frac{190781} {100297} , ~~  \frac{466872} {245443} ,
$$
$$ \frac{1124525} {591183} , ~~  \frac{1591397} {836626} , ~~  \frac{56823420} {29873093} , ~~  \frac{58414817} {30709719} , ~~  \frac{115238237} {60582812} , ~~  \frac{173653054} {91292531} , ~~
$$
$$
 \frac{462544345} {243167874} , ~~  \frac{2023830434} {1063964027} , ~~  \frac{8557866081} {4499023982} , ~~  \frac{10581696515} {5562988009} , ~~  \frac{19139562596} {10062011991} \bigg)
$$
From these denominators $q_n(\BBB)$   ~we can calculate the quantities $L(n)$:
\bee
L(n)= \left(q_n(\BBB)\right)^{1/n}, ~~~n=1, 2, \ldots , 23
\eee
The obtained values of $L(n)$  for $n\geq 7$ are presented in the third column
of the Table 1 and are shown in the Fig.2.  These  values scatter around
the red line representing the Khinchin---L{\`e}vy's constant again suggesting
that \BB \ is irrational.

In conclusion we can say that to draw firmer statement much more digits of
the Brun's constant are needed.

\bigskip
\bigskip
{\bf Acknowledgement:} The calculation were performed using  the free package PARI/GP
 \cite{Pari} developed especially for number theoretical
purposes and which has built in a set of function to deal with continued fractions.

\begin{figure}
\vspace{-0.3cm}
\begin{center}
\includegraphics[width=\textwidth,angle=0]{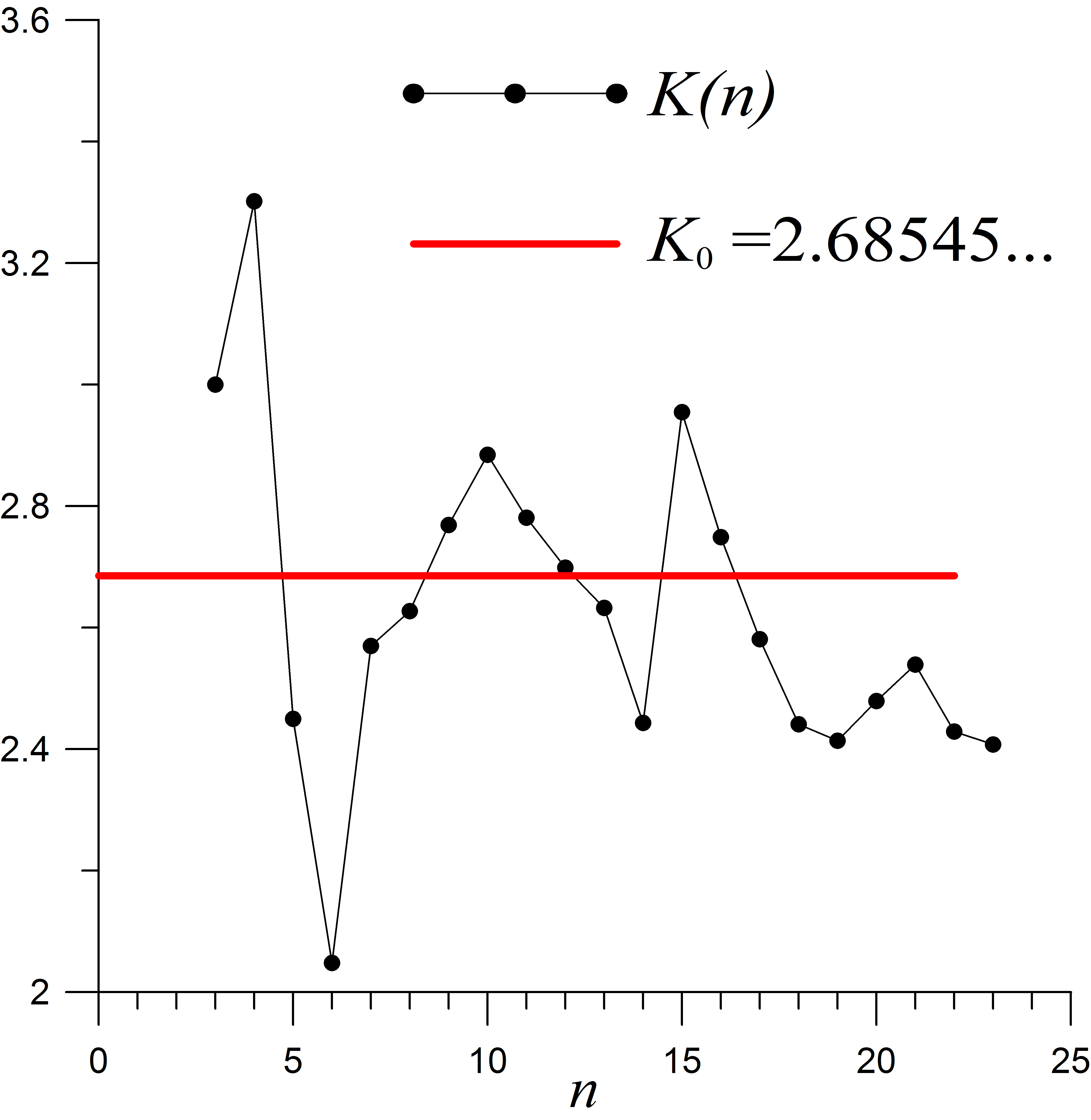} \\
\vspace{1.7cm} Fig.1 The plot of the consecutive geometrical means $K(n)$,
$n=3, 4, \ldots, 23$.
Although the number of available for the value (\ref{Sebah}) of \BB points
$(N, K(n))$
is rather moderate there are 5 sign changes of the difference $K(n)-K_0$.  \\
\end{center}
\end{figure}

\begin{figure}
\vspace{-0.3cm}
\begin{center}
\includegraphics[width=\textwidth,angle=0]{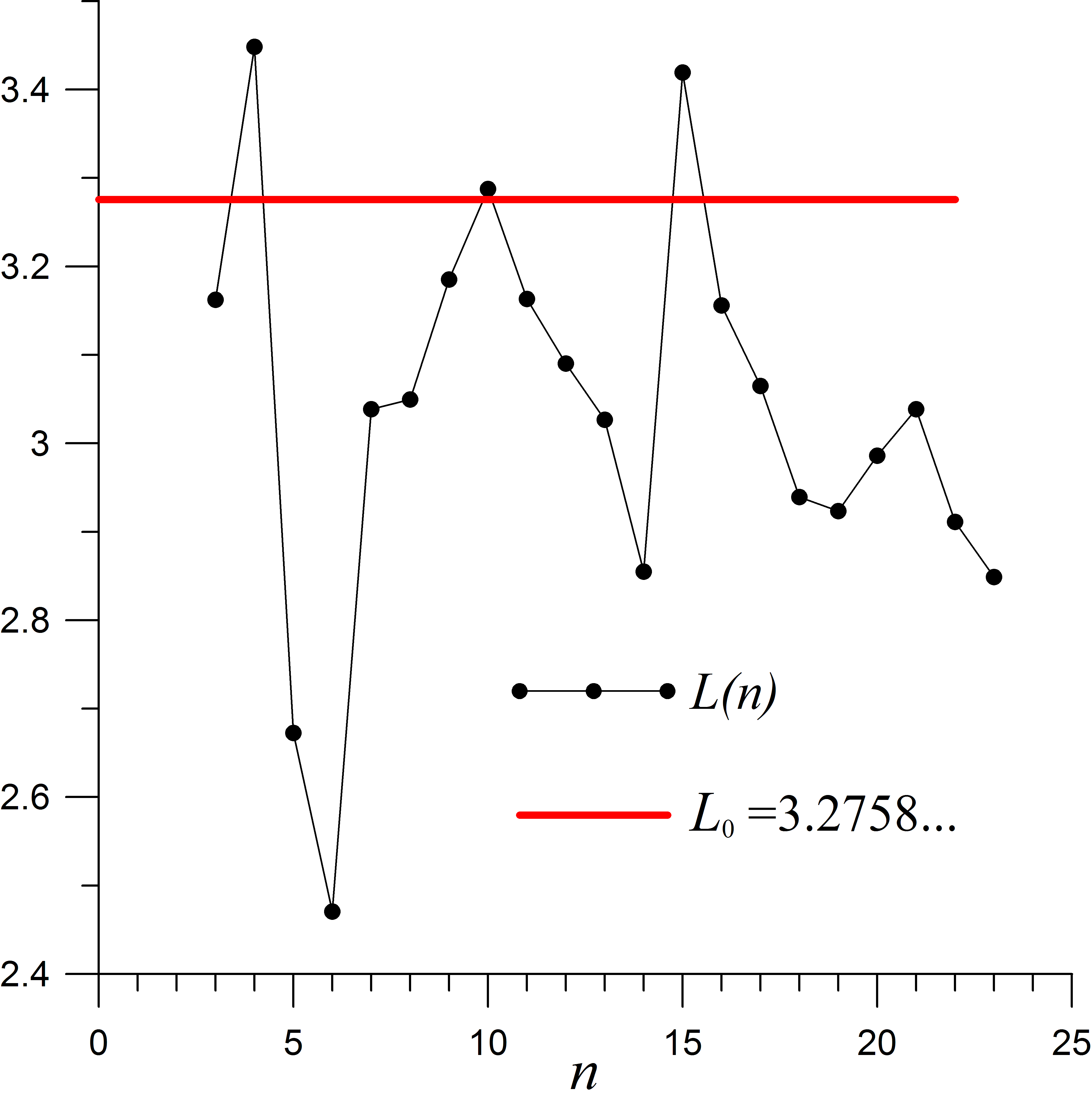} \\
\vspace{1.7cm} Fig.2 The plot of the consecutive values of  $q_n^{1/n}$,
$n=3, 4, \ldots, 23$.
Although the number of available for the value (\ref{Sebah}) of \BB points
$(n, L(n))$
is rather moderate there are 6 sign changes of the difference $L(n)-L_0$.  \\
\end{center}
\end{figure}

\end{document}